\documentclass[12pt]{article}

\title{On the Structure of Modular Categories}
\author{Michael M\"uger\thanks{Supported by NWO.} \\
Faculteit Wiskunde en Informatica, Universiteit Utrecht, Netherlands \\
and Korteweg-de Vries Institute, Amsterdam, Netherlands \\ 
email: {\tt mmueger@science.uva.nl}}

\newlength{\dinwidth}
\newlength{\dinmargin}
\usepackage{latexsym, amssymb, amsmath, theorem,tan-gle}

\def\1#1{{\bf #1}}
\def\2#1{{\cal #1}}
\def\3#1{{\sl #1}}
\def\4#1{{\tt #1}}
\def\5#1{{\sf #1}}
\def\6#1{{\mathfrak #1}}
\def\7#1{{\mathbb #1}}

\newcommand{\ba}{\begin{array}}
\newcommand{\ea}{\end{array}}
\newcommand{\bea}{\begin{eqnarray}}
\newcommand{\eea}{\end{eqnarray}}
\newcommand{\bean}{\begin{eqnarray*}}
\newcommand{\eean}{\end{eqnarray*}}
\newcommand{\nn}{\nonumber}

\newcommand{\impl}{\Rightarrow}
\newcommand{\rarr}{\rightarrow}

\newcommand{\ol}{\overline}

\newcommand{\id}{\mbox{id}}
\newcommand{\obj}{\mbox{Obj}}
\newcommand{\mcirc}{\,\circ\,}
\newcommand{\gal}{\mbox{Gal}}

\newcommand{\Hom}{\mbox{Hom}}
\newcommand{\End}{\mbox{End}}

\newcommand{\mmod}{\mbox{mod}}

\def\endexem{\hfill{$\Box$}\medskip}
\newcommand{\qed}{\ \hfill $\blacksquare$\medskip}

\newcommand{\npb}{Nucl. Phys. \1B}
\newcommand{\cmp}{Commun. Math. Phys. }

\theoremheaderfont{\scshape}
\newtheorem{defin}{Definition}[section]
\newtheorem{lemma}[defin]{Lemma}
\newtheorem{prop}[defin]{Proposition}
\newtheorem{theorem}[defin]{Theorem}
\newtheorem{coro}[defin]{Corollary}
\newtheorem{conj}[defin]{Conjecture}
\theorembodyfont{\rmfamily}
\newtheorem{rema}[defin]{Remark}

\newcommand{\bdefin}{\begin{defin}}
\newcommand{\blemma}{\begin{lemma}}
\newcommand{\bprop}{\begin{prop}}
\newcommand{\btheor}{\begin{theorem}}
\newcommand{\bcoro}{\begin{coro}}
\newcommand{\edefin}{\end{defin}}
\newcommand{\elemma}{\end{lemma}}
\newcommand{\eprop}{\end{prop}}
\newcommand{\etheor}{\end{theorem}}
\newcommand{\ecoro}{\end{coro}}
\newcommand{\bconj}{\begin{conj}}
\newcommand{\econj}{\end{conj}}
\newcommand{\brem}{\begin{rema}}
\newcommand{\erem}{\endexem\end{rema}}

\newcommand{\prf}{{\it Proof. }}

\def\mobj#1{\raise .4\unitlens\hbox{\put(0,0){$#1$}}}

\def\mychi{\raise 2pt\hbox{$\chi$}}

\begin{document}
\maketitle\noindent

\numberwithin{equation}{section}

\abstract{For a braided tensor category $\2C$ and a subcategory $\2K$ there is a notion of
centralizer $C_\2C(\2K)$, which is a full tensor subcategory of $\2C$. A pre-modular
tensor category \cite{brug1} is known to be modular in the sense of Turaev iff the center 
$\2Z_2(\2C)\equiv C_\2C(\2C)$ (not to be confused with the center $\2Z_1$ of a tensor 
category, related to the quantum double) is trivial, i.e.\ consists only of multiples of
the tensor unit, and $\dim\2C\ne 0$. Here $\dim\2C=\sum_i d(X_i)^2$, the $X_i$ being the 
simple objects. 

We prove several structural properties of modular categories. Our main technical tool is
the following double centralizer theorem. Let $\2C$ be a modular category and $\2K$ a
full tensor subcategory closed w.r.t.\ direct sums, subobjects and duals.
Then $C_\2C(C_\2C(\2K))=\2K$ and $\dim\2K\cdot\dim C_\2C(\2K)=\dim\2C$.

We give several applications. (1) If $\2C$ is modular and $\2K$ is a full modular
subcategory, then also $\2L=C_\2C(\2K)$ is modular and $\2C$ is equivalent as a ribbon 
category to the direct product: $\2C\simeq\2K\boxtimes\2L$. Thus every modular category
factorizes (non-uniquely, in general) into prime ones. We study the prime factorizations
of the categories $D(G)-\mmod$, where $G$ is a finite abelian group.
(2) If $\2C$ is a modular $*$-category and $\2K$ is a full tensor subcategory then 
$\dim\2C\ge\dim\2K\cdot\dim \2Z_2(\2K)$. We give examples where the bound is attained and
conjecture that every pre-modular $\2K$ can be embedded fully into a modular category
$\2C$ with $\dim\2C=\dim\2K\cdot\dim \2Z_2(\2K)$. (3) For every finite group $G$ there is a
braided tensor $*$-category $\2C$ such that $\2Z_2(\2C)\simeq\mbox{Rep}\,G$ (thus
$Gal(\2C)\cong G$) and the modular closure/modularization \cite{mue06,brug1}
$\ol{\ol{\2C}}$ is non-trivial.}


\section{Introduction}
Braided tensor categories \cite{js1} play a central r\^{o}le in the representation theory
of quantum groups \cite{cp,ka}, of Kac-Moody algebras \cite{kl} and of quantum field
theories \cite{frs1,ms2}. They also serve as input data for the construction of invariants
of knots, links and 3-manifolds \cite{t,ka}. In both areas -- representation theory and low 
dimensional topology -- a particular subclass of braided tensor categories is distinguished,
that of modular categories. First formalized in \cite{t0}, they are semisimple rigid
ribbon categories that have finitely many isomorphism classes of simple objects and
satisfy a non-degeneracy condition. Modular categories derive their name from the fact
that they define \cite{t} a (projective) finite dimensional representation of $SL(2,\7Z)$. 
At first sight mysterious, this modular representation is best understood in topological
terms. Viz., a modular category gives rise, for every closed oriented surface $M$, to a
finite dimensional projective representation of the mapping class group of $M$. The
modular representation associated with a modular category is then just the representation
of the mapping class group of the torus. In the important special case of $*$-categories
(or unitary categories) arising in quantum field theory this modular representation had
been known earlier and rigorously studied in \cite{khr1}.

The non-degeneracy condition mentioned above amounts to non-degeneracy of a certain
matrix, which is just the collection of invariants of the Hopf link for the possible
labelings of the two components. Also the meaning of the non-degeneracy condition in the
construction of the 3-manifold invariant is quite transparent. Yet, it is clearly 
desirable to have a more intrinsic understanding in purely algebraic terms. In the special
case of unitary categories, it has long been known \cite{khr1} that the modularity
condition is equivalent to the absence of `degenerate' objects. The latter property has a
very satisfactory interpretation in terms of triviality of the center $\2Z_2(\2C)$ of the
braided category. This center is a canonical full symmetric subcategory of $\2C$ and must
not be confused with another notion of `center' which is defined for any -- not
necessarily braided -- tensor category. We denote the latter by $\2Z_1(\2C)$. 
In the more general situation of a pre-modular category \cite{brug1}, the equivalence
between modularity and triviality of the center has been proven only recently
\cite{brug3,bebl}. (In \cite[Corollary 7.11]{mue10} this proof is obtained as a
byproduct.)

Since symmetric tensor categories are precisely those braided tensor categories which
coincide with their center, modular categories may be seen as braided tensor categories
diametrically opposed to the symmetric ones: Modular categories are related to symmetric 
tensor categories like groups with trivial center to abelian groups, or like factors to
commutative von Neumann algebras. Now, under some additional conditions, symmetric tensor
categories are just representation categories of groups \cite{del,dr6}, thus they should
be considered as very basic algebraic objects. Our point of view is that modular
categories merit to be perceived similarly and to be subjected to detailed scrutiny and, as
far as feasible, classification. So far, very little was known in the way of a general
theory, our Theorem \ref{t-factor}, according to which every modular category is
(equivalent to) a finite direct product of prime ones, apparently being one of the first 
structural results.

On the other hand, many different constructions of modular categories are known, which we
briefly review. There is a large class of constructions which go under the heading of
quantum doubles, all of which yield modular categories under suitable assumptions. For the
definition of the quantum double of finite dimensional Hopf algebras (in particular, group
algebras) and of tensor categories (the `center' $\2Z_1$ mentioned above) we refer
to \cite{ka}. Proofs of modularity were given in \cite{ac1} for quantum doubles of finite
groups and in \cite{ac2} for the twisted versions of \cite{dpr}, and in \cite{eg} and
\cite[Appendix]{mue10} for semisimple cosemisimple Hopf algebras. These cases are subsumed
by quantum groupoids and by semisimple spherical categories, for which modularity of the
quantum double was proven in \cite{ntv} and \cite{mue10}, respectively. (Actually, by
results of Hayashi \cite{hay} and Ostrik \cite{ostr}, the categories considered in
\cite{mue10} are always representation categories of quantum groupoids. Thus the proof of
the main theorem of \cite{mue10} could also be deduced from the results of \cite{ntv} with
some additional effort. Conceptually, however, the direct proof seems more satisfactory.)
Quantum groups at roots of unity give rise to modular categories if one considers
appropriate quotients of their representation categories \cite{tw}. These categories are
in fact $*$-categories \cite{w} (or unitary categories \cite{t}). (The categories obtained
in this way and from quantum doubles, respectively, are not completely unrelated since the
universal $R$-matrix of a deformed enveloping algebra $A$ is computed by expressing $A$ as
a quotient of a quantum double. On the categorical side, every modular category $\2C$ is a
full tensor subcategory of its own quantum double $\2Z_1(\2C)$, cf.\ \cite{mue10}. This is
remarkable insofar as the definition of the latter does not refer to the braiding of
$\2C$.) For some of the above modular categories a beautiful purely combinatorial
construction is known \cite{bebl}. Furthermore, there is an operation \cite{brug1,mue06}
which, heuristically, amounts to dividing a pre-modular category $\2C$ by its center and
which can be interpreted as a Galois completion \cite{mue06,brug2}. This procedure is
applicable whenever the objects of the center $\2Z_2(\2C)$ have positive integer
dimensions (automatic for $*$-categories) and trivial twists, and it
yields a modular category that is non-trivial whenever $\2C$ is not symmetric. Finally,
for a suitable class of rational chiral conformal field theories, axiomatized using
operator algebras, one can prove \cite{klm} that the category of representations is a
modular $*$-category. (In WZW- and orbifold models the representation categories are those
of \cite{w} and of \cite{dpr}, respectively.) For a review of these results we refer to 
\cite{mue08}, where also some of the results of the present paper were announced. 

In the next section we briefly review some of the formalism of modular categories,
restricting ourselves to the facts that are needed in this paper. 
With the exception of Lemma \ref{charfct} the results are well known, but we
emphasize the r\^{o}le of centralizers in braided categories. In Section 3 we prove our
main technical result, a double centralizer theorem in modular categories. Section 4 gives
several applications, the most important of which is that every modular subcategory of a
modular category is a direct factor. This implies that every modular category factorizes
as a finite direct product of prime ones. Here, a modular category $\2C$ is prime if
every full modular subcategory is either trivial or equivalent to $\2C$.
In Section 5 we give some preliminary results about full embeddings of pre-modular
categories into modular categories and conclude with a remark about the Galois theory for
braided tensor categories developed in \cite{mue06,brug2}.


\section{Preliminaries on Modular Categories}
\subsection{Notation}
We assume known the standard definitions of (braided, symmetric) tensor categories, cf.\
\cite{cwm, js1}. All categories in this paper are supposed small and all tensor categories
strict. (A tensor category
is strict if the tensor product satisfies associativity $X\otimes(Y\otimes Z)=(X\otimes
Y)\otimes Z$ `on the nose' and the unit object $\11$ satisfies $X\otimes\11=\11\otimes
X=X\ \forall X$. By the coherence theorems, every tensor category is equivalent to a
strict one.) Our notation is fairly standard, except that we omit the $\otimes$-symbol for
the product of objects: $XY\equiv X\otimes Y$. Unit objects and unit morphisms in tensor
categories are denoted by $\11$ and $\id_X\in\End(X)$, respectively. Our categories will
be $\7F$-linear over a field $\7F$. We use capital letters $X,Y,\ldots$ to denote both
objects and isomorphism classes of simple objects. What is meant should be obvious from
the context. As usual, we denote $N_{XY}^Z=\dim\Hom(Z,XY)$. A $*$-category \cite{glr,lro}
or unitary category \cite{t} is a $\7C$-linear category equipped with an antilinear 
involutive and contravariant endofunctor $*$ that leaves the objects fixed and such that 
$s^*\circ s=0$ implies $s=0$.

For the definitions of (left) rigid and ribbon categories see \cite{js1,t}. Every left
rigid ribbon category is spherical \cite{bw}, i.e.\ every object $X$ has a two-sided dual
$\ol{X}$ and left and right traces coincide. (Conversely, every spherical category
with braiding automatically has a compatible ribbon structure \cite{y}.) We denote the
twist by $\{\theta_X, X\in\2C\}$. To every simple object $X$ in a ribbon category we
assign $\omega_X\in\7F$ by $\omega_X\id_X=\theta_X$. One has $\omega_X=\omega_{\ol{X}}$
for all simple $X$, and in a $*$-category, $|\omega_X|=1$ for all $X$. A pre-modular
category \cite{brug1} is a semisimple $\7F$-linear rigid ribbon category with finitely
many isomorphism classes of simple objects and tensor unit satisfying
$\End\,\11\cong\7F$. 

A subcategory $\2S\subset\2C$ is full iff 
$\Hom_\2S(X,Y)=\Hom_\2C(X,Y) \ \ \forall X,Y\in\2S$, thus it is determined by
$\obj\,\2S$. A subcategory $\2S$ is replete iff $X\in\2S$ implies $Y\in\2S$ for all
$Y\in\2C$ isomorphic to $X$. 
By a semisimple tensor subcategory of a semisimple spherical category $\2C$ we mean
a full subcategory which is stable w.r.t.\ direct sums, subobjects (thus in particular
replete) and duals. If $\Gamma$ denotes the set of isomorphism classes of simple
objects of $\2C$, the semisimple subcategories are in one-to-one correspondence with the
subsets $\Gamma'\subset\Gamma$ which are closed under duals and satisfy
$X,Y\in\Gamma', N_{XY}^Z\ne 0\ \impl\ Z\in\Gamma'$.


\subsection{Monodromies in braided tensor categories}
\bdefin The monodromy of two objects $X,Y$ in a tensor category with braiding $c$ is
defined by
\[ c_M(X,Y)=c(Y,X)\circ c(X,Y) \ \in \End(XY). \]
\edefin

\bdefin \label{Y-m}
For a braided spherical category $\2C$ over $\7F$ and $X,Y\in\2C$ define $S(X,Y)\in\7F$ by
\[ S(X,Y)\id_\11=Tr_{XY}(c_M(X,Y))\ = \  \begin{tangle}
\hcoev\Step\hcoev\\
\id\step\xx\step\id\\
\id\hstep\mobj{X}\hstep\xx\mobj{Y}\step\id\\
\hev\step\step\hev
\end{tangle}
\]
\edefin

\brem 1. $S(X,Y)$ depends only on the isomorphism classes $[X],[Y]$.

2. Note that we did not assume $X,Y$ to be simple. 
\erem

\blemma \label{lem1}
Let $\2C$ be a $\7F$-linear semisimple rigid ribbon category with $\End\,\11\cong\7F$.
The following identities hold:
\begin{itemize}
\item[(i)] $\displaystyle S(UX,Y)=\frac{1}{d(Y)}\,S(U,Y)
   S(X,Y)\quad\forall Y\ \mbox{simple}, \forall U,X$.
\item[(ii)] $\displaystyle X\cong\bigoplus_i X_i,\ Y\cong
  \bigoplus_j Y_j\quad \Longrightarrow \quad S(X,Y)=\sum_{i,j} 
  S(X_i,Y_j)$.
\item[(iii)] $\displaystyle \frac{1}{d(X)} S(X,Y)S(X,Z)=
   \sum_W N_{YZ}^{W}\,S(X,W)\quad
  \forall X,Y,Z$ simple.
\end{itemize}
\elemma
\prf The first claim follows from $Tr_{UX}=Tr_X\circ (Tr_U\otimes\id)$ together with the 
consequence $(Tr_U\otimes\id)(c_M(U,Y))=d(Y)^{-1}S(U,Y)\id_Y$ of simplicity of $Y$. The
second claim is immediate by cyclic invariance of the trace. Part (iii) follows by
applying (i) and (ii) to $S(UV,Y)$. 
\qed 

The next result is valid without the restriction to unitary (i.e.\ $*$-) categories, cf.\
\cite{brug3, bebl}. The proof uses a general result on handle slides \cite{ab}, which we
do not wish to enlarge upon. Rather than referring the result away completely, we give a
simple argument for unitary categories which does not require $\2C$ to be finite.

\bprop \label{p-lem2}
Let $\2C$ be a $\7F$-linear semisimple rigid ribbon category with $\End\,\11\cong\7F$ 
(equivalently, a tensor $*$-category with direct sums, subobjects and conjugates).
Let $X,Y$ be simple. Then $S(X,Y)=d(X)d(Y)$ iff $c_M(X,Y)=\id_{XY}$.  
\eprop
\prf The `if' statement is obvious from the definition of $S$. Thus assume
$S(X,Y)=d(X)d(Y)$. Using the well known equations \cite{t,khr1}
\bean d(X)d(Y) &=& \sum_Z N_{XY}^Z\,d(Z), \\
   S(X,Y) &=& \sum_Z N_{XY}^Z\frac{\omega_Z}{\omega_X\omega_Y} d(Z), 
\eean
this implies
\[ \sum_Z N_{XY}^Z d(Z)\, \frac{\omega_Z}{\omega_X\omega_Y}=\sum_Z N_{XY}^Z d(Z). \]
Restricting the summation to those $Z$ for which $N_{XY}^Z\ne 0$ we have
\begin{equation} {\sum_Z}' N_{XY}^Z d(Z)
   \left(1-\frac{\omega_Z}{\omega_X\omega_Y}\right)=0.  
\label{x77}\end{equation}
Since the $\omega's$ have absolute value one, we have 
$\mbox{Re} (1-\omega_Z/\omega_X\omega_Y)\ge 0$, with equality iff 
$\omega_Z/\omega_X\omega_Y=1$. In view of $N_{XY}^Z d(Z)>0$ in (\ref{x77}) we conclude 
that $\omega_Z=\omega_X\omega_Y$ whenever
$N_{XY}^Z>0$. Let $Z$ be simple and consider $s: XY\rarr Z,\ t: Z\rarr XY$. Then with the
ribbon condition
$\theta_{XY}=\theta_X\otimes\theta_Y\circ c_M(X,Y)=c_M(X,Y)\circ \theta_X\otimes\theta_Y$
we have
\bean c_M(X,Y)\mcirc t\circ s &=& \theta_{XY}\circ(\theta_X\otimes\theta_Y)^{-1}\mcirc
   t\circ s  \ = \ (\omega_X\omega_Y)^{-1}\, \theta_{XY}\mcirc t\circ s  \\
  &=& (\omega_X\omega_Y)^{-1}\, t\mcirc\theta_{Z}\mcirc s
   \ =\ \frac{\omega_Z}{\omega_X\omega_Y}\, t\circ s. 
\eean
Since $\End(XY)$ is unital and spanned by morphisms $t\circ s$ as above, the fact that 
$\omega_Z=\omega_X\omega_Y$ for all $Z$ contained in $XY$ implies $c_M(X,Y)=\id_{XY}$. 
\qed

\bdefin Let $\2C$ be a braided tensor category and $\2K$ a set of objects in $\2C$,
equivalently, a full subcategory of $\2C$. Then we define the centralizer $C_\2C(\2K)$ of
$\2K$ in $\2C$ (or relative commutant $\2C\cap\2K'$) as the full subcategory defined by
\[ \obj \, C_\2C(\2K)=\{X\in\2C\ | \  c_M(X,Y)=\id_{XY}\ \ \forall Y\in\2K \}. \]
\edefin

\brem 1. In \cite[Subsection 5.2]{mue06} the subcategory $C_\2C(\2K)$ was called 
$\2C_\2K$. In subfactor theory, a related notion appears under the name `permutant' in
\cite{ocn3}.

2. If there is no danger of confusion concerning the ambient category $\2C$, we will
occasionally write $\2K'$ instead of $C_\2C(\2K)$.
\erem

If $\2K_1,\2K_2$ are full subcategories of $\2C$, by $\2K_1\vee \2K_2$ we denote the
smallest replete full subcategory of $\2C$ containing $\2K_1$ and $\2K_2$ and stable under
tensor products, direct sums and retractions.  

\blemma \label{l-centr}
For $\2C, \2K$ as above, $C_\2C(\2K)$ is replete and monoidal. If $\2C$ is
a semisimple category with direct sums and subobjects, then the same holds for 
$C_\2C(\2K)$. If $\2C$ has duals for all objects then also $C_\2C(\2K)$ has
duals. If $\2K_1,\2K_2\subset\2C$ then 
$C_\2C(\2K_1\vee\2K_2)=C_\2C(\2K_1)\cap C_\2C(\2K_2)$.
\elemma
\prf Repleteness of $C_\2C(\2K)$ follows from naturality of the braiding $c$. It
is easy to see that $c_M(X_i,Y)=\id_{X_iY}$ for $i=1,2$ implies
$c_M(X_1X_2,Y)=\id_{X_1X_2Y}$, thus $C_\2C(\2K)$ is closed under tensor products. In a 
semisimple category with $X=\oplus_{i\in I}X_i$ one has $ c_M(X,Y)=\id_{XY}$ iff 
$c_M(X_i,Y)=\id_{X_iY}$ for all $i\in I$. This implies that $C_\2C(\2K)$ is closed under 
direct sums and subobjects. The statement concerning duals follows by the same
argument as in the proof of \cite[Proposition 2.7]{mue06}. As to the last claim, the
inclusion $C_\2C(\2K_1\vee\2K_2)\subset C_\2C(\2K_1)\cap C_\2C(\2K_2)$ is obvious. If $X$
is in $C_\2C(\2K_1)$ and $C_\2C(\2K_2)$ then it also has trivial monodromy $c_M$ with all 
tensor products of objects in $\2K_1$ and $\2K_2$, as well as direct sums and retracts of
such. 
\qed

\bdefin \label{center}
The center of a braided tensor category $\2C$ is
\[ \2Z_2(\2C)=C_\2C(\2C). \]
We say a semisimple braided tensor category has trivial center, if every object of
$\2Z_2(\2C)$ is a direct sum of copies of the tensor unit $\11$ or if, equivalently, every
simple object in $\2Z_2(\2C)$ is isomorphic to $\11$.
\edefin

\brem \label{r-cen}
1. Clearly, a braided tensor category is symmetric iff it coincides with its 
center. 

2. The objects of the center have previously been called {\it degenerate}
\cite{khr1,mue06}, {\it transparent} \cite{brug1} and {\it pseudotrivial} \cite{saw}. Yet,
calling them {\it central} seems the most natural terminology, since the above definition
is the correct analogue for braided tensor categories of the center of a monoid, as can be
seen appealing to the theory of $n$-categories. 

3. Proposition \ref{p-lem2} now has the interpretation that the simple object $X$ is
central iff $S(X,Y)=d(X)d(Y)$ for all simple $Y\in\2C$.
\erem


\subsection{Finite dimensional categories} 
In our considerations so far we did not make finiteness assumptions on $\2C$. From now on
we will work with categories which have finite dimension in the following sense.

\bdefin \label{dimcat}
Let $\2C$ be a semisimple $\7F$-linear spherical category with $\End(\11)\cong\7F$.
If the set $\Gamma$ of isomorphism classes of simple objects is finite, the {\it
dimension} of $\2C$ is defined by 
\[ \dim\2C=\sum_{X\in\Gamma} d(X)^2, \]
otherwise it is $\infty$.
\edefin

\brem 1. The sum over the squared dimensions appeared in \cite{t}, where its square roots
are called ranks of the category. In subfactor theory \cite{ocn1,ocn2} this number is
called the `global index'. The designation `dimension' for this number has also been used
in \cite{bw}. It is vindicated by the fact that $\dim\mbox{Rep}(H)=\dim H$ for a finite 
dimensional semisimple Hopf algebra, in particular group algebras. That $\dim\2C$ is the
correct generalization of $\dim H$ it is corroborated by its behavior under various
constructions like the quantum double, where $\dim\2Z_1(\2C)=(\dim\2C)^2$ \cite{mue10}. 

2. The dimension of a semisimple $\7F$-linear category can be defined unambiguously
whenever the category has two-sided duals, cf.\ \cite[Subsect.\ 2.2]{mue09}. A
sovereign/spherical or $*$-structure is not needed.
\erem

If $\2K$ is a subcategory of $\2C$, let $\mychi_\2K$ be the characteristic function of
$\obj\,\2K$. Viz., $\mychi_\2K(X)=1$ if $X\in\2K$ and $\mychi_\2K(X)=0$ otherwise. 

\blemma \label{charfct}
Let $\2C$ be a pre-modular category and $\2K$ a semisimple tensor subcategory. 
Then for all $X\in\2C$ we have
\begin{equation} \label{x0}
  \sum_{Y\in\2K} d(Y) S(X,Y)=d(X) \dim\2K \ \mychi_{C_\2C(\2K)}(X). 
\end{equation} 
\elemma
\prf If $X\in C_\2C(\2K)$ then $S(X,Y)=d(X)d(Y)$ for all $Y\in\2K$, and (\ref{x0}) follows
immediately. Thus it remains to show that the left hand side of (\ref{x0}) vanishes if 
$ c_M(X,Y)\ne\id_{XY}$ for some $Y\in\2K$. To this purpose, consider Lemma \ref{lem1} (c)
with $Y,Z\in\2K$. Since $\2K$ is a sub-tensor category, the summation runs only over
isomorphism classes of simple $W\in\2K$.
Multiplication with $d(Y)$ and summation over $Y\in\2K$ yields
\[ \frac{S(X,Z)}{d(X)}\sum_{Y\in\2K} d(Y) 
  S(X,Y)= \sum_{Y,W\in\2K} d(Y) 
   N_{YZ}^{W}\,S(X,W). \]
Now, $\sum_{Y} d(Y) N_{YZ}^{W}=\sum_{Y} d(Y) N_{\ol{W}Z}^{\ol{Y}}=d(W)d(Z)$
and we obtain
\[ \left(\sum_{Y\in\2K} d(Y) S(X,Y)\right)
   [S(X,Z)-d(X)d(Z)]=0 
   \quad\quad\forall X\in\2C,\ Z\in\2K.  \]
If now $X\not\in C_\2C(\2K)$ then there exists a $Z\in\2K$ which has non-trivial 
monodromy with $X$, and by Proposition \ref{p-lem2} we have $S(X,Z)\neq d(X)d(Z)$. 
Thus the expression in the big brackets vanishes and we are done. \qed

For $\2K=\2C$, the lemma reduces to a known result, and together with Remark \ref{r-cen}
it implies the following corollary, cf.\ e.g.\ \cite{bebl}. But for our purposes the
generalization to arbitrary tensor subcategories $\2K$ will be essential. 

\bcoro \label{p-central}
Let $\2C$ be a pre-modular category with $\dim\2C\ne 0$ and let $X$ be a simple object.
Then the following are equivalent: 
\begin{itemize}
\item[(i)] $X$ is central.
\item[(ii)] $S(X,Y)=d(X)d(Y)$ for all simple $Y\in\2C$.
\item[(iii)] $\displaystyle \sum_Y S(X,Y)d(Y)\ne 0$.
\end{itemize}
\ecoro

\blemma \label{p-scalar}
Let $\2C$ be pre-modular and $Y,Z$ simple. Then
\[ \sum_X S(X,Y)S(X,Z) = \dim\2C \sum_{W\in \2Z_2(\2C)} N_{YZ}^{W}\,d(W). \]
If $\2Z_2(\2C)$ is trivial then $S^2=\dim\2C\,C$, where 
$C=(C_{XY}), C_{XY}=\delta_{X,\ol{Y}}$.
\elemma
\prf Multiplying (iii) of Lemma \ref{lem1} with $d(X)$ and summing over $X$ we obtain 
\bean \sum_X S(X,Y)S(X,Z) &=&  \sum_{X,W} N_{YZ}^{W}\,S(X,W) d(X) \nn\\
  &=& \dim\2C  \sum_{W\in \2Z_2(\2C)} N_{YZ}^{W}\,d(W),
\eean
where we have used Lemma \ref{charfct} with $\2K=\2C$. The last claim follows from 
$N_{XY}^0=\delta_{X,\ol{Y}}$.
\qed

\bcoro \label{c-modular} 
Let $\2C$ be pre-modular with $\dim\2C\ne 0$. Then the following are equivalent:
\begin{itemize}
\item[(i)] The center $\2Z_2(\2C)$ is trivial.
\item[(ii)] The matrix $(S(X,Y))$, indexed by isomorphism classes of simple objects, is
invertible. 
\end{itemize}
\ecoro
\prf The implication (ii)$\impl$(i) is obvious. Conversely, if $\2Z_2(\2C)$ is trivial,
the lemma gives $S^2=\dim\2C\,C$. (ii) then follows from invertibility of $C$ and
$\dim\2C\ne 0$.
\qed

\brem 1. Alternatively, the statement of the corollary can be obtained directly from a
handle slide argument, cf.\ \cite{bebl, mue10}.

2. If $\2C$ is a $*$-category, one easily shows $\ol{S(X,Y)}=S(X,\ol{Y})=S(\ol{X},Y)$.
In the case where $\2Z_2(\2C)$ is trivial, Lemma \ref{p-scalar} implies that the columns
(or rows) of $S$ are mutually orthogonal, thus $\tilde{S}=(\dim\2C)^{-1/2}S$ is unitary.
Without modularity one can still show that for simple $X,Y$ the columns $S_X,S_Y$ (or
rows) of the $S$-matrix are either orthogonal or parallel, cf.\ \cite{khr1}. If
$\ol{\ol{\2C}}$ is the modularization and $F(X)$ the image of $X$ in $\ol{\ol{\2C}}$,
t.f.a.e.: (i) $S_X\|S_Y$, (ii) there exists $Z\in\2Z_2(\2C)$ such that
$\Hom(X,ZY)\ne\{0\}$, (iii) $\Hom_{\ol{\ol{\2C}}}(F(X),F(Y))\ne\{0\}$, (iv) $F(X),F(Y)$
have the same simple summands, cf.\ \cite{mue06}.
\erem


\section{Centralizers in Modular Categories}
\bdefin \cite{t0} A modular category is a pre-modular category satisfying the (equivalent)
conditions of Corollary \ref{c-modular}. \edefin

We are now in a position to prove our first main result.

\btheor \label{dct}
Let $\2C$ be a modular category and let $\2K\subset\2C$ be a semisimple tensor
subcategory. Then we have
\begin{itemize}
\item[(i)] $\displaystyle C_\2C(C_\2C(\2K))=\2K$.
\item[(ii)] $\displaystyle \dim \2K \,\cdot\, \dim\,C_\2C(\2K)=\dim\2C$.
\end{itemize}
(We also simply write $\2K''=\2K$ and $\dim\2K\cdot\dim\2K'=\dim\2C$.)
\etheor
\prf We apply Lemma \ref{charfct} to compute the characteristic function of
$C_\2C((C_\2C(\2K))$:
\bean \mychi_{C_\2C(C_\2C(\2K))}(X) &=& \frac{1}{d(X)\dim\,C_\2C(\2K)}
   \sum_{Z\in C_\2C(\2K)} S(X,Z) d(Z) \\
  &=& \frac{1}{d(X)\dim\,C_\2C(\2K)} \sum_{Z\in\2C} \mychi_{ C_\2C(\2K)}(Z) 
   \,S(X,Z) d(Z). 
\eean
We use the lemma once again to compute $\mychi_{ C_\2C(\2K)}(Z)$, obtaining
\bean \mychi_{C_\2C((C_\2C(\2K))}(X) &=& \frac{1}{d(X)\dim\,C_\2C(\2K)}
   \sum_{Z\in\2C} S(X,Z) d(Z) \ \frac{1}{d(Z)\dim\2K} \sum_{U\in\2K}
    S(Z,U) d(U) \\
  &=& \frac{1}{d(X)\dim\,C_\2C(\2K) \dim\2K} \sum_{U\in\2K} d(U) \sum_{Z\in\2C}
   S(X,Z) S(Z,U).
\eean
The summation over $Z\in\2C$ can be performed using Lemma \ref{p-scalar}, and 
since the center of $\2C$ is trivial, by Corollary \ref{c-modular} we have
\[ \sum_{Z\in\2C} S(X,Z) S(Z,U)= \dim\2C\,\delta_{[X],\ol{[U]}}. \] 
Using $d(U)=d(\ol{U})$ and the fact that $\2K$ is closed w.r.t.\ duals we obtain
\bea \mychi_{C_\2C(C_\2C(\2K))}(X) &=& \frac{\dim\2C}{d(X)\dim\2K\dim\,C_\2C(\2K)} \
  \sum_{U\in\2K} d(U)\delta_{[X],\ol{[U]}} \nn \\
  &=& \frac{\dim\2C}{\dim\2K\dim\,C_\2C(\2K)} \ \mychi_{\2K}(X). \label{e-dct}
\eea
Since the tensor unit $\11$ is contained in any tensor subcategory we have  
$\mychi_\2K(\11)=\mychi_{\2K''}(\11)=1$. Thus for $X=\11$, (\ref{e-dct}) proves
claim (ii), and plugging this back into (\ref{e-dct}), claim (i) ensues.
\qed

\brem 1. Using (i), one easily verifies the following: If $\2K$ is any subcategory closed 
w.r.t.\ duals then $C_\2C(C_\2C(\2K))$ is the semisimple tensor subcategory generated
by $\2K$, i.e.\ the completion w.r.t.\ direct sums and subobjects.

2. In a subfactor context, the double centralizer property $\2K''=\2K$ was stated
by A. Ocneanu \cite{ocn3} without published proof. Subfactor analogues of both (i) and  
(ii) were proved by Izumi \cite{iz2} using considerable machinery. By contrast, the above 
proof uses only well known properties modular categories.
\erem

\bcoro \label{cor1}
Let $\2C$ be a modular category and let $\2K\subset\2C$ be a semisimple tensor
subcategory. Then
\ecoro
\[ \2Z_2( C_\2C(\2K))=\2Z_2(\2K). \]
\prf We compute
\bean \lefteqn{
  \2Z_2(C_\2C(\2K))=C_{C_\2C(\2K)}(C_\2C(\2K))=(\2C\cap\2K')\cap(\2C\cap\2K')' } \\
  &&= (\2C\cap(\2C\cap\2K')')\cap\2K'=C_{C_\2C(C_\2C(\2K))}(\2K)=C_\2K(\2K)=\2Z_2(\2K), 
\eean
where we have used $C_\2C(C_\2C(\2K))=\2K$. \qed

The two most interesting cases are where $\2K$ is modular or symmetric.

\bcoro \label{c-coro2}
Let $\2C$ be a modular category and let $\2K\subset\2C$ be a semisimple tensor subcategory
that is modular. Then $\2L\equiv C_\2C(\2K)$ is modular, too.
\ecoro
\prf By Corollary \ref{c-modular}, modularity of $\2K$ is equivalent to triviality of 
$\2Z_2(\2K)=C_\2K(\2K)$. By Corollary \ref{cor1}, also $C_\2L(\2L)=C_\2K(\2K)$ is
trivial, thus $\2L$ is modular.
\qed

\bcoro \label{c-coro3}
Let $\2C$ be a modular category and let $\2K\subset\2C$ be a semisimple tensor subcategory
that is symmetric. Let $\2L\equiv C_\2C(\2K)$. Then $\2Z_2(\2L)=\2K$.
\ecoro
\prf Obvious: $\2Z_2(\2K)=\2K$. \qed


\section{On the Structure of Modular Categories}
\subsection{Prime factorization of modular categories}
Our first application of the double centralizer theorem in modular categories is also the
most striking one. It illustrates how different modular categories are from their opposite
extreme case, viz.\ the symmetric categories, the group duals (at least in characteristic
zero \cite{dr6,del}).

The following was proved in \cite{mue10} as Corollary 7.8. (The proof uses only two
preceding results and is independent of the rest of \cite{mue10}.) Here $\2A\boxtimes\2B$
is the completion w.r.t.\ direct sums of the product of $\2A$ and $\2B$ as $\7F$-linear 
categories. $\2A\boxtimes\2B$ has an obvious tensor structure if $\2A,\2B$ do.

\bprop \label{cor-fact}
Let $\2C$ be a braided tensor category that is $\7F$-linear semisimple with
$\End\,\11\cong\7F$ and two-sided duals. Let $\2K\subset\2C$ be a semisimple tensor
subcategory which has trivial center $\2Z_2(\2K)$. Then we have the equivalence 
\[ \2K\boxtimes C_\2C(\2K) \ \simeq \ \2K\vee C_\2C(\2K) \]
of braided tensor categories. If $\2C$ is spherical then by restriction also $\2K$ and
$C_\2C(\2K)$ are spherical, and the above equivalence is one of spherical
categories. Similarly, if $\2C$ is a $*$-category.
\eprop

The question arises naturally, whether $\2K\vee\2L=\2C$, since this would imply
$\2C\simeq\2K\boxtimes\2L$. If $\2C$ is modular, the double centralizer theorem provides
the missing step.

\btheor \label{tprod}
Let $\2C$ and $\2K$ be modular categories where $\2K$ is identified with a full (tensor) 
subcategory of $\2C$. Let $\2L= C_\2C(\2K)$. Then there is an equivalence of ribbon
categories:
\[ \2C \ \simeq \ \2K\boxtimes\2L. \]
\etheor
\prf Modularity of $\2L$ has been proved in Lemma \ref{c-coro2}. If we can show that the
full subcategories $\2K$ and $\2L$ of $\2C$ generate $\2C$, Proposition \ref{cor-fact}
provides an equivalence $\2C\simeq\2K\boxtimes\2L$ of braided tensor categories. The
equivalence is automatically an equivalence of ribbon categories since $\2K$ and $\2L$
commute. (Alternatively, one appeals to the compatibility of the spherical structures.)
For the remaining fact $\2K\vee\2L=\2C$ we give two proofs, the first of which works 
only for unitary modular categories. 

{\it Unitary categories.} By Proposition \ref{cor-fact}, we have
$\dim\,\2K\vee\2L=\dim\2K\cdot\dim\2L$, which by Theorem \ref{dct} coincides with
$\dim\2C$. Since $\2K\vee\2L$ is a full subcategory of $\2C$ and the numbers
$d(X)^2\in\7R$ are non-negative, the equality
\[ \dim\2C=\sum_{X\in\2C} d(X)^2=\sum_{X\in\2K\vee\2L} d(X)^2=\dim\,\2K\vee\2L \]
implies that all simple objects of $\2C$ are contained in $\2K\vee\2L$ and therefore
$\2C=\2K\vee\2L\simeq\2K\boxtimes\2L$, as desired.

{\it General case.} The preceding argument to the effect that $\2K\vee\2L$ exhausts
$\2C$ does not work if we are not dealing with $*$-categories, since $\dim\2A=\dim\2B$
does not imply that a replete full inclusion $\2A\subset\2B$ is an identity. Yet, we can 
argue as follows:
\bean \lefteqn{ \2K\vee\2L=C_\2C(C_\2C(\2K\vee\2L))=C_\2C(C_\2C(\2K)\cap C_\2C(\2L)) 
  =C_\2C(C_\2C(\2K)\cap C_\2C(C_\2C(\2K))) } \\
 && =C_\2C(C_\2C(\2K)\cap \2K) =C_\2C(\2C\cap\2K'\cap\2K)
    =C_\2C(C_\2K(\2K))=C_\2C(\2Z_2(\2K))=\2C. 
\eean
Here we used Lemma \ref{l-centr} and the fact that $\2Z_2(\2K)$ is trivial.
\qed

\brem 1. That $\2L$ is modular can be derived alternatively from the easy fact that a
direct product $\2C\simeq\2K\boxtimes\2L$ is modular iff both $\2K$ and $\2L$ are modular.

2. An interesting special case is provided by the quantum double of a modular
category $\2K$. If we put $\2C=\2Z_1(\2K)$, there exists a braided monoidal embedding
functor $I: \2K\hookrightarrow\2C$. The above theorem implies 
$\2C\simeq I(\2K)\boxtimes C_{\2C}(I(\2K))$. In this situation, $C_\2C(I(\2K))$ can be
computed explicitly and is given by $\tilde{I}(\tilde{\2K})$, where $\tilde{\2K}$
coincides with $\2K$ as a tensor category, but has the opposite braiding, and $\tilde{I}$
is a braided monoidal embedding. Thus $\2C\simeq\2K\boxtimes\tilde{\2K}$, which is the
statement of \cite[Theorem 7.9]{mue10}. Even if $\2K$ is not modular, $I(\2K)$ 
and $\tilde{I}(\tilde{\2K})$ are each others centralizers in $\2C$, and we have
$I(\2K)\cap\tilde{I}(\tilde{\2K})=I(\2Z_2(\2K))$. See \cite[Section 7]{mue10} for the
details.

3. Finally, it is interesting to note the analogy with the following result from
classical non-commutative algebra. Let $A\subset B$ be an inclusion of matrix algebras.
Then (i) $C_B(A)$ is a matrix algebra, (ii) $C_B(C_B(A))=A$, and (iii) 
$B\cong A\otimes C_B(A)$. These results also hold in the infinite dimensional case if
$A,B$ are type I factors (i.e.\ isomorphic to the algebras of bounded linear operators on
some Hilbert spaces).
\erem

\bdefin A modular category $\2C$ is prime if every semisimple modular subcategory is
equivalent either to $\2C$ or the trivial modular category $\mbox{Vect}_\7F$.
\edefin

Now we can state our second main result.
\btheor \label{t-factor}
Every modular category is equivalent to a finite direct product of prime ones. 
\etheor
\prf Obvious consequence of Theorem \ref{tprod} since modular categories have finitely
many (equivalence classes of) simple objects and proper replete full subcategories have
strictly fewer simple objects.
\qed

\brem 1. A. Brugui\`{e}res \cite{brug3} first observed that the above factorization in
general is non-unique. Yet, uniqueness does hold if every simple object except the unit
has dimension unequal to one. In the next subsection we will find many examples for the
non-uniqueness. 

2. For the classification of modular categories, Theorem \ref{t-factor} has the obvious
consequence that it is sufficient to consider prime ones. As in other respects, modular
categories are better behaved than finite groups since there are no non-trivial exact
sequences to be considered. 
\erem


\subsection{Quantum doubles of finite abelian groups}
Quantum doubles of finite groups provide a large and relatively easy-to-analyze class of
modular categories, cf., e.g., \cite{ac1}. As an exploratory step towards a classification
of modular ($*$-)categories it is natural to find the prime factorizations of the categories
$D(G)-\mmod$ ($\equiv D(\7CG)-\mmod$). If $G$ is a direct product of subgroups $K,L$ then
it is easy to see that $D(G)\cong D(K)\otimes D(L)$ as quasitriangular Hopf algebra, thus
$D(G)-\mmod\simeq D(K)-\mmod\boxtimes D(L)-\mmod$. Therefore, in order for $D(G)-\mmod$ to
be prime, $G$ must be prime (not a direct product of non-trivial subgroups). This
condition is, however, not sufficient, as is shown by the complete analysis of the abelian
case given below.

If $G$ is abelian, the simple objects of $\2C=D(G)-\mmod$ have dimension one and are
invertible, the isomorphism classes of simple objects and their tensor product being given
by the abelian group $\Gamma(\2C)=G\times\hat{G}$. There is an obvious one-to-one
correspondence between replete full tensor subcategories $\2K\subset\2C$ and subgroups
$K\subset\Gamma(\2C)$. Apart from this we only need the fact that
$S((g,\chi),(h,\sigma))=\langle\sigma,g\rangle\langle\chi,h\rangle$.

It is well known that every finite abelian group is isomorphic to a direct product of
cyclic groups of prime power order,
\[ G\cong \7Z/p_1^{n_1}\7Z\,\times \ldots\times\, \7Z/p_k^{n_k}\7Z, \]
where the pairs $(p_i,n_i)$ are unique up to permutation. By the above remark we can
therefore restrict ourselves to a consideration of $G=\7Z/p^n\7Z$. 

\btheor 
Let $p$ be prime, $G=\7Z/p^n\7Z$ and $\2C=D(G)-\mmod$.
\begin{itemize}
\item[(i)] If $p=2$ then $\2C$ is prime.
\item[(ii)] If $p$ is odd, there is a one-to-one correspondence between isomorphisms
$\alpha: G\rarr\hat{G}$ and replete full modular subcategories $\2K_\alpha\subset\2C$
given by $K_\alpha=\{(g,\alpha(g)),\ g\in G\}$. The categories $\2K_\alpha$ are prime,
and $C_\2C(\2C_\alpha)=\2C_{\ol{\alpha}}$, where $\ol{\alpha}(\cdot)=\alpha(\cdot)^{-1}$.
The prime factorizations of $\2C$ are thus given by 
$\2C\simeq\2C_\alpha\boxtimes\2C_{\ol{\alpha}}, \ \alpha\in\mbox{Isom}(G,\hat{G})$.
\end{itemize}
\etheor
\prf Let $\2K\subset\2C$ be a replete full tensor subcategory corresponding to the subgroup
$K\subset\Gamma$. If $\2K$ is modular then, by Theorem \ref{tprod}, $\Gamma\cong K\times
L$, where $L$ corresponds to the centralizer $\2L= C_\2C(\2K)$. By the uniqueness result
for the factorization of finite abelian groups, a (non-trivial) factorization
$\Gamma\cong K\times L$ of $\Gamma\cong(\7Z/p^n\7Z)^2$ is possible only if 
$K\cong L\cong\7Z/p^n\7Z$. Identifying $G$ and $\hat{G}$ with $\7Z/p^n\7Z$, the $S$-matrix
is given by
\[ S((a,b),(c,d))=e^{\frac{2\pi i}{p^n}(ad+bc)}. \]
Furthermore, every cyclic subgroup $K\subset\Gamma$ of order $p^n$ is of the form
$\{(ja,jb), j\in\7Z/p^n\7Z\}$, where $a,b$ are not both multiples of $p$. 
$\2K$ is modular iff $j\ne 0$ implies $(ja,jb)\not\in\2Z_2(\2K)$, which by Proposition
\ref{p-central} is equivalent to 
\[ \sum_{k=0}^{p^n-1} e^{\frac{2\pi i}{p^n}2jkab} = 0.\]
In view of $\sum_{k=0}^{N-1} a^k=(a^N-1)/(a-1)$ for $a\ne 1$ and $=N$ otherwise, this is
the case iff $e^{\frac{2\pi i}{p^n}2jab}\ne 1$ and $e^{\frac{2\pi i}{p^n}2jabp^n}=1$. 
The latter condition is always satisfied, and the former leads to
\begin{equation} \label{e-xxx} 
  (ja,jb)\in\2Z_2(\2K) \quad\Leftrightarrow\quad 2jab\in p^n\7Z. \end{equation}
If $p=2$ then (\ref{e-xxx}) is satisfied by $j=p^{n-1}$ irrespective of $a,b$. Thus
$(ja,jb)\in\2Z_2(\2K)$ and $\2K$ is not modular, proving (i). From now on let $p$ be odd. If
$p\mid a$ or $p\mid b$ then again $j=p^{n-1}$ solves (\ref{e-xxx}) and $\2K$ is not
modular. If $p\nmid a$ and $p\nmid b$ then (\ref{e-xxx}) is satisfied iff 
$j\equiv 0 (\mbox{mod}\ p^n)$, which implies modularity of $\2K$. Rephrasing this in an
invariant fashion leads to the first part of statement (ii). Proper subgroups of
$K_\alpha$ have order smaller than $p^n$ and therefore cannot give rise to modular
subcategories by the argument at the beginning of the proof. Thus the
$\2K_\alpha,\alpha\in\mbox{Isom}(G,\hat{G})$ are prime, and they exhaust the prime factors
of $\2C$. With $\alpha,\alpha'\in\mbox{Isom}(G,\hat{G})$, Proposition \ref{p-lem2} finally 
implies that $\2K_\alpha'=\2K_{\alpha'}$ iff 
\[ S((g,\alpha(g)),(h,\alpha'(h)))=\langle\alpha(g),h\rangle\langle\alpha'(h),g\rangle =1
   \quad\forall g,h\in G. \]
For $G=\7Z/p^n\7Z$ and $\alpha\in\mbox{Isom}(G,\hat{G})$ one easily shows
$\langle\alpha(g),h\rangle=\langle\alpha(h),g\rangle\ \forall g,h$, which evidently
entails $\2K_\alpha'=\2K_{\ol{\alpha}}$.
\qed

\brem 1. In view of $\#\mbox{Isom}(G,\hat{G})=p^n-p^{n-1}$, the theorem nicely illustrates
to which extent the factorization of modular categories into primes can be non-unique.

2. The table
\begin{center}
\begin{tabular}{lcc}
$G$ & $G$ simple? & $D(G)-\mmod$ prime? \\ \hline
$\7Z/2\7Z$ & Yes & Yes \\
$\7Z/p\7Z, \ p\ne 2$ & Yes & No \\
$\7Z/2^n\7Z, \ n\ge 2$ & No & Yes \\
$\7Z/p^n\7Z, \ p\ne 2, n\ge 2$ & No & No 
\end{tabular}
\end{center}

\noindent
makes clear that already for an abelian group $G$ there is no relation between simplicity
of $G$ and primality of $D(G)-\mmod$.
\erem


\section{Further Applications}
\subsection{Modular extensions of braided categories} \label{embed}
The theory of Galois extensions of braided tensor categories developed in
\cite{brug1,mue06} provides a means to construct a modular category $\ol{\ol{\2C}}$ from a
given pre-modular category $\2C$. In the language of \cite{mue06}, this modular closure is
given by $\ol{\ol{\2C}}=\2C\rtimes \2Z_2(\2C)$, i.e.\ by adding morphisms which turn the
objects in $\2Z_2(\2C)$ into multiples of the tensor unit $\11$. The dimension of the
modular closure is given by 
\begin{equation} \label{x1}
    \dim\,\ol{\ol{\2C}}=\frac{\dim\2C}{\dim \2Z_2(\2C)},
\end{equation}
cf. \cite{brug1,mue15}. Thus $\ol{\ol{\2C}}$ is, in a sense, a quotient of $\2C$ by
$\2Z_2(\2C)$, in fact it is trivial iff $\2C$ is symmetric.

On the other hand one may wish to find a modular category related to $\2C$ without
modifying the latter. More precisely, the problem is to find a modular category
$\hat{\2C}$, into which $\2C$ embeds as a full subcategory. If we restrict ourselves to 
$*$-categories, for which the dimensions take values in $\7R_+$, we obtain a lower bound
on the dimension of $\hat{\2C}$ as an immediate corollary of Theorem \ref{dct}. 

\bprop \label{cor2}
Let $\2C$ be a unitary modular category and $\2K\subset\2C$ a semisimple tensor
subcategory. Then 
\begin{equation} \dim\2C\ge\dim\2K\cdot\dim \2Z_2(\2K). \label{bd1}\end{equation}
Equality holds iff $C_\2C(\2K)=C_\2K(\2K)=\2Z_2(\2K)$.
\eprop
\prf The obvious inclusion (of replete full tensor subcategories) 
$C_\2C(\2K)\supset C_\2K(\2K)=\2Z_2(\2K)$ implies
$\dim\,C_\2C(\2K)\ge\dim\,\2Z_2(\2K)$ and therefore
$\dim\2C=\dim\2K\cdot\dim\,C_\2C(\2K)\ge\dim\2K\cdot\dim\2Z_2(\2K)$. Equality in
(\ref{bd1}) is equivalent to $\dim\,C_\2C(\2K)=\dim\,\2Z_2(\2K)$ and thus to
$C_\2C(\2K)=\2Z_2(\2K)$.
\qed

It is natural to expect that this bound is optimal.

\bconj
Let $\2C$ be a unitary pre-modular category. Then there exists a unitary modular category 
$\hat{\2C}$ and a full and faithful tensor functor $I: \2C\rarr\hat{\2C}$ such that
\[ \dim\hat{\2C}=\dim\2C\cdot\,\2Z_2(\2C). \] 
Such a category $\hat{\2C}$ is called a minimal modular extension of $\2C$.
\econj

\brem 1. It is instructive to compare this equation with (\ref{x1}). In fact, the
appearance of the factor $\dim\2Z_2(\2C)$ in the conjecture and in (\ref{x1}) is not
accidental. This is evident from the orbifold construction in conformal field theory,
where the Galois closure $\ol{\ol{\2C}}$ and the minimal modular extension $\hat{\2C}$
of a pre-modular category $\2C$ both appear naturally, cf.\ \cite{mue08}.

2. The conjecture makes sense also without the unitarity assumption. But without the
latter it is not clear in which sense $\hat{\2C}$ is minimal.

3. If $\2C$ is symmetric, i.e.\ $\2C=\2Z_2(\2C)$, (\ref{bd1}) reduces to
$\dim\hat{\2C}\ge(\dim\2C)^2$. If $\omega_X=1$ for all simple $X$, by \cite{dr6} we have a
finite group $G$ such that $\2C\simeq G-\mmod$ and $\dim\2C=|G|$. Since the modular
category $D(G)-\mmod$ satisfies $\dim D(G)-\mmod=|G|^2$ and contains $G-\mmod$ as a full
subcategory we see that pre-modular $*$-categories that are symmetric with trivial twists
in fact admit a minimal modular extension. 

4. The example $\2C=G-\mmod$ also shows that $\hat{\2C}$, if it exists, will in general
not be unique. The modular categories $D^\omega(G)-\mmod$ \cite{dpr}, where 
$\omega\in Z^3(G,\7T)$, are minimal modular extensions of $G-\mmod$ which are inequivalent
for different $[\omega]\in H^3(G,\7T)$.

5. In a subfactor context, the existence (and `essential uniqueness'!) of such a $\2C$ is
claimed by Ocneanu \cite{ocn3}. To the best of our knowledge a proof has not appeared. 
\erem

Despite the fact that we do not know how to construct a minimal modular extension for a
given pre-modular category, using our results \cite{mue10} on the categorical quantum
double we can construct many examples of unitary pre-modular categories which do possess a
minimal modular extension. 

Let $\2C$ be a pre-modular category with $\dim\2C\ne\infty$ and let $\2Z_1(\2C)$ be
is quantum double. 
It is known \cite{mue10} that $\2C$ is modular and satisfies $\dim\2Z_1(\2C)=(\dim\2C)^2$.
There are fully faithful braided tensor functors $I: \2C\rarr\2Z_1(\2C)$ and 
$\tilde{I}: \tilde{\2C}\rarr\2Z_1(\2C)$, where $\tilde{\2C}$ equals $\2C$ as a tensor
category, the braiding given by $\tilde{c}(X,Y)=c(Y,X)^{-1}$. Let 
\[ \2E=I(\2C)\vee\tilde{I}(\tilde{\2C}) \]
be the replete full subcategory of $\2Z_1(\2C)$ generated by $I(\2C)$ and 
$\tilde{I}(\tilde{\2C})$.

\bprop The subcategory $\2E\subset\2Z_1(\2C)$ satisfies
\begin{itemize}
\item[(i)] $\displaystyle \2Z_2(\2E)=I(\2Z_2(\2C))=\tilde{I}(\2Z_2(\tilde{\2C}))$.
\item[(ii)] $\displaystyle\dim\2E=\frac{(\dim\2C)^2}{\dim \2Z_2(\2C)}$.
\end{itemize}
\eprop
\prf We compute the centralizer of $\2E$ in $\2Z_1(\2C)$:
\bean C_{\2Z_1(\2C)}(\2E) &=& C_{\2Z_1(\2C)}(I(\2C)\vee\tilde{I}(\tilde{\2C})) =
  C_{\2Z_1(\2C)}(I(\2C))\cap C_{\2Z_1(\2C)}(\tilde{I}(\tilde{\2C})) \\ 
  &=& \tilde{I}(\tilde{\2C})\cap I(\2C)= I(\2Z_2(\2C)). \eean
(We have used Lemma \ref{l-centr} and the observations from \cite[Section 4]{mue10} that
$I(\2C)$ and $\tilde{I}(\tilde{\2C})$ are each other's centralizer in $\2Z_1(\2C)$ and 
satisfy $I(\2C)\cap\tilde{I}(\tilde{\2C})=I(\2Z_2(\2C))$.) This computation obviously
implies $\2Z_2(\2E)=\2E\cap\2E'=I(\2Z_2(\2C))$. Furthermore, we have
$\dim\,C_{\2Z_1(\2C)}(\2E)=\dim\,\2Z_2(\2C)$, and using part (b) of the double centralizer
theorem we obtain
\[ \dim \2E=\frac{(\dim\2C)^2}{\dim\2Z_2(\2C)}. \]
\qed

\brem Fact (ii) is intuitively clear: If $\2Z_2(\2C)$ is trivial then $\2E$
is equivalent to the direct product of $I(\2C)\cong\2C$ and 
$\tilde{I}(\tilde{\2C})\cong\tilde{\2C}$, cf.\ \cite[Theorem 7.10]{mue10}.
Thus for non-trivial $\2Z_2(\2C)$, $\2E$ should be a product of $I(\2C)$ and
$\tilde{I}(\tilde{\2C})$ amalgamated over the common subcategory $I(\2Z_2(\2C))$. 
\erem

\bcoro $\2Z_1(\2C)\supset\2E$ is a minimal modular extension of $\2E$. \ecoro
\prf By construction, $\2Z_1(\2C)$ is modular and contains $\2E$ as replete full braided
subcategory. In view of the preceding results and the known facts on the double we have 
\[ \dim\2Z_1(\2C)=(\dim\2C)^2=\frac{(\dim\2C)^2}{\dim\2Z_2(\2C)}\dim\2Z_2(\2C)
   =\dim\2E\cdot\dim\2Z_2(\2E),  \] 
thus the bound in Proposition \ref{cor2} is attained. \qed


\subsection{The inverse problem for Galois extensions of premodular categories}
In this subsection we will solve the inverse problem of the Galois theory for
braided tensor categories which was developed in \cite{mue06}, see also \cite{brug2}.
To begin with, given an arbitrary compact group $G$, we have 
$\gal(\mbox{Rep}(G))\cong G$. Thus every compact group is the Galois group of some
braided tensor category. This is, however, not very interesting since symmetric tensor
categories have trivial modular closure. We will therefore show that for every finite
group there is a finite dimensional unitary pre-modular category with $\gal(\2C)\cong G$
and which is not symmetric, thus has non-trivial modular closure. To this purpose we will
construct modular categories $\2C$ containing $\2S=\mbox{Rep}(G)$ as a full subcategory. 
Defining $\2D=C_\2C(\2S)$, Corollary \ref{c-coro3} implies $\2Z_2(\2D)=\2S$, thus we
have $\gal(\2D)\cong G$ as desired. Non-triviality of 
$\2D\rtimes\2Z_2(\2D)=\2D\rtimes\2S$ is equivalent to $\2D\supsetneq\2S$, thus
$\dim\2C>(\dim\2S)^2$.

We exhibit to ways of constructing such a $\2C$, both of which involve the center
$\2Z_1(\2C)$ or quantum double of a tensor category $\2C$, which was already referred to
in Subsection \ref{embed}. See \cite{ka} for a nice discussion of the center construction
and \cite{mue10} for additional results which we will need.

The first construction is: $\2C=\2Z_1(\2Z_1(\2S))$. Since $\2S$ is braided, $\2Z_1(\2S)$ and
$\2Z_1\circ \2Z_1(\2S)$ contain $\2S$ as a full subcategory. By the results of \cite{mue10},
$\2C$ is modular and $\dim\2C=(\dim\2S)^4$. In particular, $\dim\2C>(\dim\2S)^2$ as
required. Since we obtain 
\[ \dim(C_\2C(\2S)\rtimes\2S)=(\dim\2S)^2=|G|^2, \]
it is natural to conjecture that the modular closure of $\2Z_1=C_\2C(\2S)$ is 
equivalent to $\2Z_1(\2S)\cong D(G)-\mbox{mod}$.

The other procedure is as follows. Let
\[ \11\longrightarrow H\longrightarrow K\longrightarrow G\longrightarrow\11 \]
be an exact sequence of finite groups with $H\ne\{e\}$.
With $\2C=\2Z_1(\mbox{Rep}(K))$, $\2C$ is
modular and $\dim\2C=|K|^2$. Furthermore, $\2S=\mbox{Rep}(G)$ is contained as a
full subcategory in $\mbox{Rep}(K)$ and thus in $\2Z_1(\mbox{Rep}(K))$. We obtain
\[ \dim(C_\2C(\2S)\rtimes\2S)=\frac{|K|^2}{|G|^2}=|H|^2. \]
We therefore expect that $\2D\rtimes\2S=C_\2C(\2S)\rtimes\2S$ is equivalent to
$\2Z_1(\mbox{Rep}(H))\cong D(H)-\mbox{mod}$. 

Thus we have given two different proofs of the following
\btheor Let $G$ be a finite group. Then there is a unitary pre-modular category $\2C$ 
such that $\mbox{Gal}(\2C)\cong G$ and such that $\2C$ is not symmetric, thus having
a non-trivial modular closure. \etheor

\vspace{1cm}\noindent
{\it Acknowledgments.} I thank A. Brugui\`{e}res and M. Izumi for useful conversations.

\end{document}